\documentclass[leqno]{amsart}

\usepackage{amsmath}
\usepackage{amsfonts}
\usepackage{amssymb}
\usepackage{eucal}

\begin{document}

\title{On pseudo-Riemannian manifolds \\ 
with recurrent concircular curvature tensor}

\author{Karina Olszak and Zbigniew Olszak}

\date{\today}

\begin{abstract}
It is proved that every concircularly recurrent manifold must be necessarily a recurrent manifold with the same recurrence form.
\end{abstract}



\address{Institute of Mathematics and Computer Science, Wroc{\l}aw University of Technology, Wybrze\.ze Wyspia\'nskiego 27, 50-370 Wroc{\l}aw, Poland}

\email{Karina.Olszak@pwr.wroc.pl, Zbigniew.Olszak@pwr.wroc.pl}

\maketitle

\begin{center}
\bf 1. Introduction
\end{center}

\medskip
For a pseudo-Riemannian manifold $(M,g)$, by $T,U,V,W,X,Y,Z$ will be denoted arbitrary smooth vector fields on $M$, and for the Riemann curvature operator $\mathcal R$ and the Riemann curvature $(0,4)$-tensor $R$, we assume the following conventions 
\begin{eqnarray*}
  & \mathcal R(X,Y)=[\nabla_X,\nabla_Y]-\nabla_{[X,Y]} 
                   = \nabla^2_{X,Y} - \nabla^2_{Y,X}, &\\[+3pt]
  & R(W,X,Y,Z) = g(\mathcal R(W,X)Y,Z), &
\end{eqnarray*}
where $\nabla$ indicates the covariant derivative with respect to the Levi-Civita connection, and $\nabla^2$ is the second covariant derivative, 
$$
  \nabla^2_{X,Y} = \nabla_X\nabla_Y - \nabla_{\nabla_XY}.
$$
Pseudo-Riemannian manifolds are assumed to be connected. 

A pseudo-Riemannian manifold $(M,g)$ is said to be recurrent \cite{RWW,Wa} if its Riemann curvature operator $\mathcal R$ is recurrent, that is, $\mathcal R$ is non-zero and its covariant derivative $\nabla\mathcal R$ satisfies the condition 
\begin{equation}
\label{rec}
  \nabla\mathcal R = \lambda\otimes\mathcal R
\end{equation}
for a certain 1-form $\lambda$ (the recurrence form). 

For a pseudo-Riemannian manifold $(M,g)$, the concircular curvature tensor field $\mathcal C$ is defined as
\begin{equation}
\label{conc}
  \mathcal C = \mathcal R - \frac{r}{n(n-1)}\mathcal G,
\end{equation}
where $n=\dim M$, $r$ is the scalar curvature and $\mathcal G$ is the curvature like operator defined as
$$
  \mathcal G(X,Y)Z = g(Y,Z)X - g(X,Z)Y.
$$

The tensor $\mathcal C$ is an invariant of the concircular transformations which have many important geometric and algebraic applications; see \cite{Y,V,KV,KP}, etc. For our purpose, we recall only two facts: (1) when $\dim M=2$, then $\mathcal C=0$ and such a manifold realizes the condition $\nabla\mathcal R=\lambda\otimes\mathcal R$ with $\lambda=\nabla(\mathop{\rm ln}|r|)$ at each point at which $\mathcal R\not=0$; (2) when $\dim M\geqslant 3$, $\mathcal C=0$ if and only if the manifold is of constant sectional curvature. 

A pseudo-Riemannian manifold $(M,g)$ is said to be concircularly recurrent if its concircular curvature tensor $\mathcal C$ is recurrent, that is, $\mathcal C$ is non-zero and its covariant derivative $\nabla\mathcal C$ satisfies the condition 
\begin{equation}
\label{concrec}
  \nabla\mathcal C = \lambda\otimes\mathcal C
\end{equation}
for a certain 1-form $\lambda$. For a concircularly recurrent manifold, $n=\dim M\geqslant3$ and $\mathcal C\not=0$ at every point of $M$. 

It is obvious that a recurrent manifold is concircularly flat ($\mathcal C=0$) or concircularly recurrent with the same recurrence form. The purpose of the presented paper is to prove the following theorem: 

\medskip
{\sc Theorem}. {\it Every concircularly recurrent manifold is necessarily a recurrent manifold with the same recurrence form.}

\medskip
The above theorem seems to be very important since concircularly recurrent manifolds were studied by many authors, see \cite{ADMY,BKT,DGa,DGu,D,MM,MR,SI}, etc. In view of our theorem, many results and proofs occured in some of the listed papers can be simplified, sometimes radically. 

\medskip
{\sc Hint}. {\it With the recurrence notion used in the present paper, we follow Y.-C.\ Wong \cite{Wo1,Wo2}; cf.\ also \cite{KN}. Due to Theorem 3.8 from \cite{Wo1}, on a connected differentiable manifold endowed with an affine connection, any recurrent tensor field has no zeros.} 

\bigskip\bigskip
\begin{center}
\bf 2. Proof of the theorem
\end{center}

\medskip
For a concirculary recurrent manifold, as a consequence of (\ref{conc}) and (\ref{concrec}), we claim that the Riemann curvature operator satisfies the following condition 
\begin{equation}
\label{grec}
  \nabla\mathcal R = \lambda\otimes\mathcal R + \mu\otimes\mathcal G,
\end{equation}
where the 1-form $\lambda$ is the same as in (\ref{concrec}) and the 1-form $\mu$ is given by 
\begin{equation}
\label{mu}
  \mu = \frac{1}{n(n-1)}(dr - r\lambda). 
\end{equation}
Conversely, if a pseudo-Riemannian manifold satisfies (\ref{grec}) with certain 1-forms $\lambda$ and $\mu$, then $\mu$ must be of the form (\ref{mu}), and (\ref{concrec}) must be realized. The equivalence of (\ref{concrec}) and (\ref{grec}) has already been noticed in \cite{ADMY}. It is worth to remark that if $\mu=0$, then obviously (\ref{rec}) holds so that the concircular recurrence reduces to the recurrence. 

Before we start with the proof, it will be useful to recall the famous Walker identity (see \cite[Lemma 1]{Wa}) stating that for any pseudo-Riemannian manifold it holds
\begin{eqnarray*}
  & (\nabla^2_{U,V}R-\nabla^2_{V,U}R)(W,X,Y,Z) 
     + (\nabla^2_{W,X}R-\nabla^2_{X,W}R)(Y,Z,U,V) &\\
  & \null+ (\nabla^2_{Y,Z}R-\nabla^2_{Z,Y}R)(U,V,W,X) = 0. &
\end{eqnarray*}
Rewrite the Walker identity in the following form, which will be more convenient for us 
\begin{equation}
\label{walk}
  (\mathcal R(U,V)R)(W,X,Y,Z) + (\mathcal R(W,X)R)(Y,Z,U,V) + (\mathcal R(Y,Z)R)(U,V,W,X) = 0.
\end{equation}

We are going to reach the assertion of the Theorem in the following three steps.

However, we would like to add that the closedness of the form $\lambda$ (see the first step below) has already been proved in \cite{ADMY}. We have included it only for the completness of the whole of the proof. 

\medskip
\noindent
{\bf Step 1}.
{\it For a concircularly recurrent manifold, the recurrence form $\lambda$ is closed.}

\begin{proof}
Let $C$ denotes the $(0,4)$-tensor related to $\mathcal C$ by
$$
  C(W,X,Y,Z) = g(\mathcal C(W,X)Y,Z). 
$$
Using (\ref{concrec}), we have $\nabla_V C = \lambda(V) C$, and next  
$$
  \nabla^2_{U,V}C=((\nabla_U\lambda)(V)+\lambda(U)\lambda(V)) C.
$$
Therefore, 
$$
  \mathcal R(U,V)C = \nabla^2_{U,V}C - \nabla^2_{V,U}C = 2d\lambda(U,V)C,
$$
where we have used the formula 
$$
  d\lambda(U,V)=\frac{1}{2}((\nabla_U\lambda)(V)-(\nabla_V\lambda)(U)). 
$$
Consequently, we obtain 
\begin{eqnarray}
\label{concwalk}
  && (\mathcal R(U,V)C)(W,X,Y,Z) + (\mathcal R(W,X)C)(Y,Z,U,V) \\
  && \quad\null + (\mathcal R(Y,Z)C)(U,V,W,X) = 2d\lambda(U,V)C(W,X,Y,Z) \nonumber\\
  && \quad\null+ 2d\lambda(W,X)C(Y,Z,U,V) + 2d\lambda(Y,Z)(Y))C(U,V,W,X),\nonumber
\end{eqnarray}
On the other hand, note that from (\ref{conc}) it follows that 
\begin{equation}
\label{conc04}
  C = R - \frac{r}{n(n-1)}\,G,
\end{equation}
where $G$ is the curvature like $(0,4)$-tensor related to $\mathcal G$ by 
$$
  G(W,X,Y,Z) = g(\mathcal G(W,X)Y,Z).
$$
Using (\ref{conc04}) and the identities $\mathcal R(U,V)r=0$, $\mathcal R(U,V)G=0$, we find 
\begin{eqnarray}
\label{twolefts}
  && (\mathcal R(U,V)C)(W,X,Y,Z) + (\mathcal R(W,X)C)(Y,Z,U,V) \\
  && \quad\null + (\mathcal R(Y,Z)C)(U,V,W,X) = (\mathcal R(U,V)R)(W,X,Y,Z) \nonumber\\
  && \quad\null + (\mathcal R(Y,Z)R)(U,V,W,X) + (\mathcal R(W,X)R)(Y,Z,U,V). \nonumber 
\end{eqnarray}
Therefore, from (\ref{walk}), (\ref{twolefts}) and (\ref{concwalk}), we have 
$$
  d\lambda(U,V)C(W,X,Y,Z) + d\lambda(W,X)C(Y,Z,U,V) + d\lambda(Y,Z)C(U,V,W,X) = 0.
$$
Since $C\not=0$, by applying the famous Walker lemma (\cite[Lemma 2]{Wa}), we obtain from the last formula 
\begin{equation}
\label{dlam}
  d\lambda = 0, 
\end{equation}
completing the proof of Step 1. 
\end{proof}

\medskip
\noindent
{\bf Step 2}.
{\it A concircularly recurrent manifold is semisymmetric.}

\begin{proof}
Using (\ref{grec}), for the first covariant derivative of the Riemann curvature $(0,4)$-tensor $R$, we find 
\begin{equation}
\label{grec2}
  \nabla_V R = \lambda(V)R + \mu(V)G, 
\end{equation}
and for the second covariant derivative, 
$$
  \nabla^2_{U,V}R = ((\nabla_U\lambda)(V)+\lambda(U)\lambda(V))R
  + ((\nabla_U\mu)(V)+\mu(U)\lambda(V))G.
$$
Hence, using also (\ref{dlam}), we obtain
\begin{equation}
\label{almsemisym}
  \mathcal R(U,V)R = \nabla^2_{U,V}R - \nabla^2_{V,U}R = 2(d\mu+\mu\wedge\lambda)(U,V)G.
\end{equation}
Applying the above identity into the Walker identity (\ref{walk}), we obtain
\begin{eqnarray*}
  & (d\mu+\mu\wedge\lambda)(U,V)G(W,X,Y,Z) 
     + (d\mu+\mu\wedge\lambda)(W,X)G(Y,Z,U,V) &\\
  & \null+ (d\mu+\mu\wedge\lambda)(Y,Z)G(U,V,W,X) = 0. &
\end{eqnarray*}
Since $G\not=0$, by applying the famous Walker lemma, we obtain
$$
  d\mu+\mu\wedge\lambda = 0.  
$$
The last relation reduces (\ref{almsemisym}) to 
\begin{equation}
\label{semisym}
  \mathcal R(U,V)R = 0,
\end{equation}
which is just the semisymmetry (cf.\ \cite{BKV}). 
\end{proof}

\noindent
\noindent
{\bf Step 3}.
{\it A concircularly recurrent manifold is recurrent.}

\begin{proof}
Since $\mathcal R(U,V)$ is a derivation of the tensor algebra on $M$ (cf.\ e.g.\ \cite{R}), we have 
\begin{eqnarray*}
  (\mathcal R(U,V)R)(W,X,Y,Z) 
  &=& \null- R(\mathcal R(U,V)W,X,Y,Z) - R(W,\mathcal R(U,V)X,Y,Z) \\
  & & \null- R(W,X,\mathcal R(U,V)Y,Z) - R(W,X,Y,\mathcal R(U,V)Z). 
\end{eqnarray*}
Hence, having the semisymmetry condition (\ref{semisym}), we obtain 
\begin{eqnarray}
\label{conseq}
  && R(\mathcal R(U,V)W,X,Y,Z) + R(W,\mathcal R(U,V)X,Y,Z) \\
  && \qquad\null+ R(W,X,\mathcal R(U,V)Y,Z) + R(W,X,Y,\mathcal R(U,V)Z) = 0. \nonumber
\end{eqnarray}
Now, differentiating the above equality covariantly, we get
\begin{eqnarray*}
  && (\nabla_TR)(\mathcal R(U,V)W,X,Y,Z) 
     + R((\nabla_T\mathcal R)(U,V)W,X,Y,Z) \\
  && \null + (\nabla_TR)(W,\mathcal R(U,V)X,Y,Z) 
     + R(W,(\nabla_T\mathcal R)(U,V)X,Y,Z) \\
  && \null + (\nabla_TR)(W,X,\mathcal R(U,V)Y,Z) 
     + R(W,X,(\nabla_T\mathcal R)(U,V)Y,Z) \\
  && \null + (\nabla_TR)(W,X,Y,\mathcal R(U,V)Z) 
     + R(W,X,Y,(\nabla_T\mathcal R)(U,V)Z) \\ 
  && \null = 0.
\end{eqnarray*}
Hence, by applying (\ref{grec}) and (\ref{grec2}), we find 
\begin{eqnarray}
\label{after14}
  && \\
  && (\lambda(T)R+\mu(T)G)(\mathcal R(U,V)W,X,Y,Z) 
     + R((\lambda(T)\mathcal R+\mu(T)\mathcal G)(U,V)W,X,Y,Z) \nonumber\\
  && \null + (\lambda(T)R+\mu(T)G)(W,\mathcal R(U,V)X,Y,Z) 
     + R(W,(\lambda(T)\mathcal R+\mu(T)\mathcal G)(U,V)X,Y,Z) \nonumber\\
  && \null + (\lambda(T)R+\mu(T)G)(W,X,\mathcal R(U,V)Y,Z) 
     + R(W,X,(\lambda(T)\mathcal R+\mu(T)\mathcal G)(U,V)Y,Z) \nonumber\\
  && \null + (\lambda(T)R+\mu(T)G)(W,X,Y,\mathcal R(U,V)Z) 
     + R(W,X,Y,(\lambda(T)\mathcal R+\mu(T)\mathcal G)(U,V)Z) \nonumber\\ 
  && \null = 0. \nonumber
\end{eqnarray}
Let us assume that $\mu\not=0$ at a certain point of $M$. At this point, using (\ref{conseq}), the equality (\ref{after14}) can be reduced to 
\begin{eqnarray*}
  && G(\mathcal R(U,V)W,X,Y,Z) + R(\mathcal G(U,V)W,X,Y,Z) \\
  && \null + G(W,\mathcal R(U,V)X,Y,Z) + R(W,\mathcal G(U,V)X,Y,Z) \\
  && \null + G(W,X,\mathcal R(U,V)Y,Z) + R(W,X,\mathcal G(U,V)Y,Z) \\
  && \null + G(W,X,Y,\mathcal R(U,V)Z) + R(W,X,Y,\mathcal G(U,V)Z) = 0.
\end{eqnarray*}
When using the definitions of the tensor $G$ and the operator $\mathcal G$, the last equality takes the following form 
\begin{eqnarray*}
  && g(V,W)R(U,X,Y,Z) - g(U,W)R(V,X,Y,Z) \\
  && \null + g(V,X)R(W,U,Y,Z) - g(U,X)R(W,V,Y,Z) \\
  && \null + g(V,Y)R(W,X,U,Z) - g(U,Y)R(W,X,V,Z) \\
  && \null + g(V,Z)R(W,X,Y,U) - g(U,Z)R(W,X,Y,V) = 0.
\end{eqnarray*}
Contracting the above with respect to the pair of arguments $V,W$ (this means that we take the trace $\mathop{\rm Trace}_g\{(V,W)\mapsto...\}$, where the dots stand for the relation to be traced), we obtain 
\begin{eqnarray*}
  & (n-2)R(U,X,Y,Z) + R(Y,X,U,Z) + R(Z,X,Y,U) & \\ 
  & \null + g(U,Y)S(X,Z) - g(U,Z)S(X,Y) = 0, &
\end{eqnarray*}
which with the help of the first Bianchi identity becomes 
\begin{equation}
\label{proj}
  (n-1)R(U,X,Y,Z) + g(U,Y)S(X,Z)-g(U,Z)S(X,Y) = 0, 
\end{equation}
$S$ being the Ricci curvature tensor. Contracting the obtained relation with respect to the pair of arguments $X,Z$, we get the Einstein condition $S = (r/n) g$, which applied to (\ref{proj}), gives us $R = (r/(n(n-1))) G$, or equivalently $C=0$, contradicting our assumption. Therefore, $\mu=0$ at every point of $M$, and consequently, the concircular recurrence reduces to the recurrence. 
\end{proof}

The authors would like to thank the reviewer for some essential remarks.

\end{document}